# Likelihood Inference in Exponential Families and Directions of Recession

## Charles J. Geyer

*School of Statistics, University of Minnesota*
*313 Ford Hall, 224 Church St. SE*
*Minneapolis, MN 55455 USA*
*e-mail:* charlie@stat.umn.edu
*url:* www.stat.umn.edu/geyer/gdor

**Abstract:** When in a full exponential family the maximum likelihood estimate (MLE) does not exist, the MLE may exist in the Barndorff-Nielsen completion of the family. We propose a practical algorithm for finding the MLE in the completion based on repeated linear programming using the R contributed package rcdd and illustrate it with two generalized linear model examples. When the MLE for the null hypothesis lies in the completion, likelihood ratio tests of model comparison are almost unchanged from the usual case. Only the degrees of freedom need to be adjusted. When the MLE lies in the completion, confidence intervals are changed much more from the usual case. The MLE of the natural parameter can be thought of as having gone to infinity in a certain direction, which we call a generic direction of recession. We propose a new one-sided confidence interval which says how close to infinity the natural parameter may be. This maps to one-sided confidence intervals for mean values showing how close to the boundary of their support they may be.

## Contents

## 1. Introduction

The problem addressed in this article is widespread. Many users of statistics have run into it, although they may not have been aware of it. The problem has been well understood for thirty years, but until now convenient software to handle the problem has not been available.

In a discrete exponential family, for example, logistic regression or categorical data analysis, it can happen that the maximum likelihood estimate (MLE) does not exist in the conventional sense. This is often not detected by software, so users may be unaware of the situation. Whether they are aware of it or not, available software provides no support for valid statistical inference in this situation. The `aster` contributed package for R (6) has a check for near-singularity of the Fisher information matrix, which indicates either nonexistence of the MLE or near collinearity of predictors, and this check has revealed a number of instances where nonexistence of the MLE arose in actual applications.

When the MLE does not exist in the conventional sense, it may exist in the Barndorff-Nielsen completion of the family (2; 3; 5, and references cited at the beginning of Chapter 2 of the latter). A practical algorithm for finding the MLE in the Barndorff-Nielsen completion of an exponential family using repeated linear programming was proposed in the author's unpublished thesis (5) and used in (10). Now we propose a different method, also using repeated linear programming with the R contributed package `rcdd` (9). We also propose methods for hypothesis tests and confidence intervals valid when the MLE does not exist in the conventional sense. Our methods only work for full exponential families but are simpler than the methods of (5) which worked for non-full convex families.

Our examples are chosen to be good examples for teaching statisticians how to use our methods. All the analysis discussed in this article is carried out in full in the accompanying technical report (7), which is produced using the R function `Sweave` so all results in the report are actually produced by the code shown therein and hence are fully reproducible by anyone who has R.

## 2. Examples

### 2.1. A Binomial Example

The binomial distribution provides the simplest example. Suppose $x$ is binomial with sample size $n$ and success probability $p$. The MLE for $p$ is $\hat{p} = x/n$. So far, so good, but three things are different about the cases $\hat{p} = 0$ and $\hat{p} = 1$. First, the natural parameter is $\theta = \text{logit}(p)$, and $\hat{\theta} = \text{logit}(\hat{p})$ does not exist when $\hat{p} = 0$ or $\hat{p} = 1$. Second, the probability distribution corresponding to the MLE is degenerate, the binomial distribution with success probability $p = 0$ or $p = 1$ being concentrated at $x = 0$ and $x = n$, respectively. Third, the elementary 95% confidence interval $\hat{p} \pm 1.96\sqrt{\hat{p}(1-\hat{p})/n}$ does not work when $\hat{p} = 0$ or $\hat{p} = 1$.

### 2.2. A Logistic Regression Example

The same kind of problems arise in multiparameter exponential families but the relevant multidimensional geometry is much harder to visualize. To introduce ideas, consider what is still a relatively simple example, which is a logistic regression. Suppose we observe a vector $y$ whose components are Bernoulli with means forming a vector $p$. The natural parameter is $\theta = \text{logit}(p)$, where logit operates componentwise $\theta_i = \text{logit}(p_i)$. Suppose we also have one covariate vector $x$ and we want to fit a quadratic model

$$\theta_i = \beta_1 + \beta_2 x_i + \beta_3 x_i^2.$$

Finally, suppose $x_i$ takes the values 1, ..., 30 and $y_i = 0$ for $x_i \le 12$ or $x_i \ge 24$ and $y_i = 1$ otherwise.

If we try to fit these data using the R function `glm`, it complains "algorithm did not converge, fitted probabilities numerically 0 or 1 occurred." In fact, the MLE for $\beta$ does not exist. Define the vector

$$\delta = (-587, 72, -2) \tag{1}$$

then we say $\delta$ is a generic direction of recession (GDOR) because there exists a $\hat{\beta}$ such that

$$\lim_{s \to \infty} l(\hat{\beta} + s\delta) = \sup_{\beta \in \mathbb{R}^3} l(\beta), \tag{2}$$

where $l$ is the log likelihood function. Although we write $\hat\beta$ here, this does not denote the MLE — the MLE does not exist — we can consider the MLE for $\beta$ to be $\hat\beta$ sent to infinity in the direction $\delta$. Neither $\hat\beta$ nor $\delta$ satisfying (2) are unique. In this example, (2) actually holds for all $\hat\beta \in \mathbb{R}^3$ when $\delta$ is given by (1).

If we consider what happens to the mean value parameter vector $p$, we find that $p(\hat\beta + s\delta) \to y$ as $s \to \infty$. Thus MLE for the mean value parameter vector does exist and is equal to the observed data vector. This MLE is strange because the distribution it goes with is degenerate, concentrated at the observed data. The MLE distribution says we can never observe data other than what we did observe; other data values occur with probability zero.

This degeneracy need not cause problems for statistical inference. The sample is not the population, and estimates are not parameters. What we need is a confidence interval, necessarily one-sided, that says how close $\beta$ is to infinity and how close $p$ is to $y$.

Fix a choice of $\hat\beta$, say $\hat\beta = 0$, and consider for all real $s$ the probability distribution having parameter $\hat\beta + s\delta$. As $s$ goes from $-\infty$ to $+\infty$ the probability of seeing the data value actually observed goes from zero to one. Find the unique $s$ that makes this probability 0.05, call it $\hat s$, then $[\hat s, \infty)$ is a 95% confidence interval for the scalar parameter $s$. Figure 1 shows the mean value parameter values corresponding to the ends of this confidence interval ($\hat s$ and $\infty$). These intervals are not the best we can do; Figure 2 shows improved intervals that require more computation.

We summarize this example, explaining which features are general and which are not. First, the GDOR notion is perfectly general. When the MLE for the natural parameter does not exist, then under conditions of Brown (3) that hold in all practical examples, there is a GDOR, and the likelihood is maximized by going to infinity in that direction. Second, the MLE for the mean value parameter for this example corresponds to a completely degenerate distribution, which is not general. Usually it will be only partially degenerate, some components of the response being fixed at an extreme value but other components being random under the MLE distribution. Third, the MLE for the natural parameter "is" any $\hat\beta$ plus infinity in the direction of the GDOR, which is not general. Usually only some $\hat\beta$ will work in (2). Fourth, Figure 1 shows most of what is important about this example, but in general there is no such figure. Usually, hypothesis tests and confidence intervals are all that can be reported.

### 2.3. A Contingency Table Example

This example is a $2 \times 2 \times \cdots \times 2$ contingency table with seven dimensions hence $2^7 = 128$ cells. The file http://www.stat.umn.edu/geyer/gdor/catrec.txt presents the data as eight vectors, seven categorical predictors $v_1$, ..., $v_7$ that specify the cells of the contingency table and one response $y$ that gives the cell counts.

To simplify the discussion, we take the cell counts to be independent Poisson random variables so we have a generalized linear model. The R statement

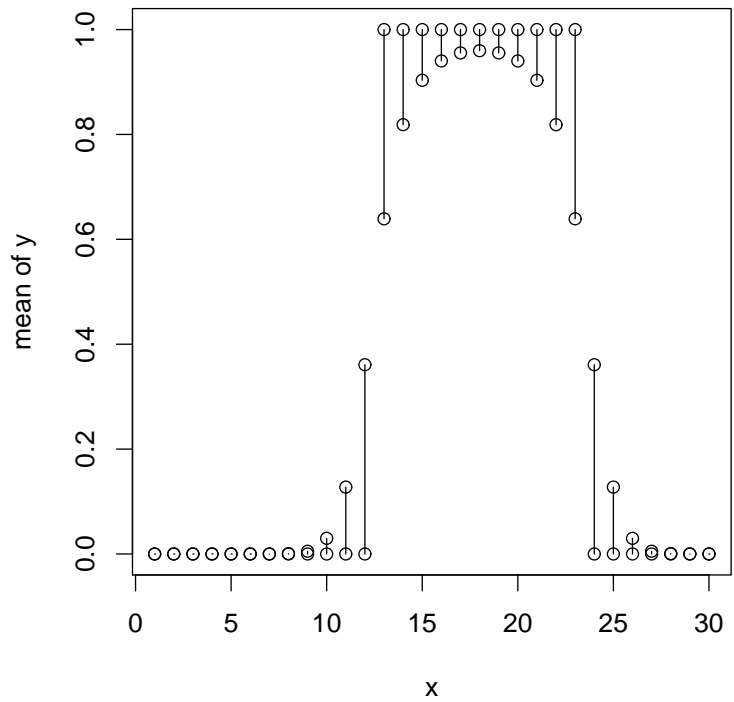

Fig 1. *Quick and dirty confidence intervals for the mean value parameter. Outer points having values zero or one are MLE for the mean value parameter, which are also the components of the observed data vector. Inner points are mean values corresponding to $\hat{\beta} + \hat{s}\delta$, where $\hat{s}$ is the lower bound for a 95% one-sided confidence interval for $s$ and where $\hat{\beta} = 0$ and $\delta$ is given by* (1). *Compare Figure 2.*

*Generic direction of recession for contingency table example. Only nonzero components shown.*

| coefficient | direction |
|---|---|
| intercept | $-1$ |
| $v_1$ | $1$ |
| $v_2$ | $1$ |
| $v_3$ | $1$ |
| $v_5$ | $1$ |
| $v_1 : v_2$ | $-1$ |
| $v_1 : v_3$ | $-1$ |
| $v_1 : v_5$ | $-1$ |
| $v_2 : v_3$ | $-1$ |
| $v_2 : v_5$ | $-1$ |
| $v_3 : v_5$ | $-1$ |
| $v_1 : v_2 : v_3$ | $1$ |
| $v_1 : v_3 : v_5$ | $1$ |
| $v_2 : v_3 : v_5$ | $1$ |

```
out3 <- glm(y ~ (v1 + v2 + v3 + v4 + v5 + v6 + v7)^3,
    family = poisson, data = dat, x = TRUE)
```

fits the model with all three-way interactions but no higher-order interactions, assuming the data have been read in as data frame `dat`.

R does not complain in fitting this model, even though it should. The correct MLE is "at infinity."

This model has 64 parameters (for a table with 128 cells). One might say this is too many parameters chasing too little data, but a test comparing this model to the model with only two-way interactions says the three-way model fits the data much better ($P \approx 10^{-17}$). Whether one likes this model or not, it should be possible for statisticians to analyze it, and we will use it for an example.

With 64 parameters we do not show the whole GDOR, only its nonzero components, which are shown in Table 1. It is hard to visualize as a 64-dimensional vector. Much easier to understand is what it does to the mean cell counts. Sixteen cells have mean zero at the MLE in the Barndorff-Nielsen completion. They are shown in Table 2. One other cell has observed data value zero but MLE mean value nonzero.

The next step in the analysis is to find a $\hat{\beta}$ such that (2) is satisfied. It has long been known that such a $\hat{\beta}$ is determined by finding the MLE in the family conditioned on certain cells being zero, in this case those in Table 2. This is easily accomplished by removing those rows from the data. The following R statements fit this model

```
dat.cond <- dat[linear, ]
out3.cond <- glm(y ~ (v1 + v2 + v3 + v4 + v5 + v6 + v7)^3,
    family = poisson, data = dat.cond)
```

where `linear` indicates the cells that have nonzero MLE mean values and `dat` is the data frame containing the original data. We do not show the results of this fit — it is shown in the accompanying technical report (7) — partly because there are 64 regression coefficients, but mainly because there is nothing new in this fit or its use to make inference about the cell mean values or the corresponding linear predictor values for the cells of the table involved (those with `linear` equal to TRUE). The R function `predict.glm` will produce valid confidence intervals for these cells using this fit.

TABLE 2
*Cells with MLE mean value zero and 95% confidence intervals (lower and upper are lower and upper confidence bounds).*

| $v_1$ | $v_2$ | $v_3$ | $v_4$ | $v_5$ | $v_6$ | $v_7$ | lower | upper |
|---|---|---|---|---|---|---|---|---|
| 0 | 0 | 0 | 0 | 0 | 0 | 0 | 0 | 0.2863 |
| 0 | 0 | 0 | 1 | 0 | 0 | 0 | 0 | 0.1408 |
| 1 | 1 | 0 | 0 | 1 | 0 | 0 | 0 | 0.2200 |
| 1 | 1 | 0 | 1 | 1 | 0 | 0 | 0 | 0.4210 |
| 0 | 0 | 0 | 0 | 0 | 1 | 0 | 0 | 0.0895 |
| 0 | 0 | 0 | 1 | 0 | 1 | 0 | 0 | 0.0938 |
| 1 | 1 | 0 | 0 | 1 | 1 | 0 | 0 | 0.1930 |
| 1 | 1 | 0 | 1 | 1 | 1 | 0 | 0 | 0.2887 |
| 0 | 0 | 0 | 0 | 0 | 0 | 1 | 0 | 0.1063 |
| 0 | 0 | 0 | 1 | 0 | 0 | 1 | 0 | 0.1141 |
| 1 | 1 | 0 | 0 | 1 | 0 | 1 | 0 | 0.0913 |
| 1 | 1 | 0 | 1 | 1 | 0 | 1 | 0 | 0.2646 |
| 0 | 0 | 0 | 0 | 0 | 1 | 1 | 0 | 0.0667 |
| 0 | 0 | 0 | 1 | 0 | 1 | 1 | 0 | 0.1548 |
| 1 | 1 | 0 | 0 | 1 | 1 | 1 | 0 | 0.1410 |
| 1 | 1 | 0 | 1 | 1 | 1 | 1 | 0 | 0.3239 |

We only remark that $\hat{\beta}$ is a little strange in that it has one nonestimable coefficient (which R indicates as being `NA`) because the model matrix is not full rank. The original model matrix was full rank, but when we remove sixteen of its rows the result is not full rank. R is smart enough to drop one or more columns (in this example just one) of the model matrix producing a new model matrix that is full rank and has the same column space as the one it was given, which is equivalent to setting the coefficients corresponding to the dropped columns to zero. Although R reports these coefficients as `NA`, we must take them to be zero.

### 2.3.1. Confidence Intervals

To obtain one-sided confidence intervals we use the same procedure used in the other example. Consider for all real $s$ the probability distribution having parameter $\hat{\beta} + s\delta$. As $s$ goes from $-\infty$ to $+\infty$ the probability of observing data having zeros in the 16 cells listed in Table 2 goes from zero to one. Find the unique $s$ that makes this probability 0.05, call it $\hat{s}$, then $[\hat{s}, \infty)$ is a 95% confidence interval for the scalar parameter $s$. Calculate the mean value parameter vector corresponding to $\hat{\beta} + \hat{s}\delta$. This gives the upper bounds for 95% confidence intervals shown in Table 2.

### 2.3.2. Hypothesis Tests

When the MLE does not exist in the conventional sense for the null hypothesis of a likelihood ratio test, the usual asymptotics do not hold and the usual test is invalid. However, the usual asymptotics may apply when the test is done conditionally, the conditioning event being the same as the one used in determining $\hat{\beta}$. In this example, if the null hypothesis is all three-way but no higher-way interactions, one merely proceeds using the data frame `data.cond` produced above to fit both the null and alternative models and do the likelihood ratio test. This idea was suggested by S. Fienberg (personal communication).

Of course, one does not have to believe the asymptotics on which this test is based. One can instead calculate a $P$-value based on a parametric bootstrap, using as the null distribution the one fitted in `out3.cond`. Since this procedure has nothing to do with the concerns of this article, being the same thing one would do whenever one distrusts asymptotics, we will say no more about it.

## 3. Theory

We redevelop the theory of Barndorff-Nielsen completion of exponential families ([2], Sections 9.3 and 9.4; [3], pp. 191–202; [5], Chapters 2 and 4) so that it is useful for calculation, particularly calculation using the R contributed package `rcdd`.

### 3.1. Exponential Families

An exponential family of distributions ([2]; [3]; [5]) is a statistical model having log likelihood

$$l(\theta) = \langle y, \theta \rangle - c(\theta), \tag{3}$$

where $y = (y_1, \ldots, y_p)$ is a vector statistic, $\theta = (\theta_1, \ldots, \theta_p)$ is a vector parameter, and

$$\langle y, \theta \rangle = \sum_{i=1}^{p} y_i \theta_i.$$

A statistic $y$ and parameter $\theta$ that give a log likelihood of this form are called *natural*. The function $c$ is called the *cumulant function* of the family.

The distribution with parameter value $\theta$ has a density with respect to the distribution with parameter value $\psi$ of the form

$$f_\theta(\omega) = e^{\langle Y(\omega), \theta - \psi \rangle - c(\theta) + c(\psi)}. \tag{4}$$

The requirement that this integrate to one determines the function $c$ up to an additive constant

$$c(\theta) = c(\psi) + \log E_\psi \big( e^{\langle Y, \theta - \psi \rangle} \big). \tag{5}$$

We take (5) to be valid for all $\theta$ in $\mathbb{R}^p$, defining $c(\theta) = \infty$ for $\theta$ such that the expectation in (5) is infinite. Define

$$\Theta = \{\, \theta \in \mathbb{R}^p : c(\theta) < \infty \,\}. \qquad (6)$$

The exponential family is *full* if its natural parameter space is (6). We shall be interested only in full families.

By convention, the cumulant function and log likelihood are defined for all $\theta \in \mathbb{R}^p$ not just at valid parameter values. We have $c(\theta) = \infty$ and $l(\theta) = -\infty$ for $\theta \notin \Theta$, so such $\theta$ can never maximize the likelihood.

### 3.2. Tangent Cone, Normal Cone, and Convex Support

The *tangent cone* of a convex set $C$ at a point $y \in C$ is

$$T_C(y) = \mathrm{cl}\{\, s(w - y) : w \in C \text{ and } s \geq 0 \,\}, \qquad (7)$$

where cl denotes the closure operation (14, Theorem 6.9). The *normal cone* of a convex set $C$ in $\mathbb{R}^p$ at a point $y \in C$ is

$$N_C(y) = \{\, \delta \in \mathbb{R}^p : \langle w - y, \delta \rangle \leq 0 \text{ for all } w \in C \,\}. \qquad (8)$$

(14, Theorem 6.9). Tangent and normal cones are polars of each other (14, Theorem 6.9 and Corollary 6.30). Each determines the other.

The *convex support* of an exponential family is the smallest closed convex set that contains the natural statistic with probability one under some distribution in the family, in which case this is true for all distributions in the family, because the distributions are mutually absolutely continuous (2, pp. 90 and 111–112).

### 3.3. Directions of Recession and Constancy

Directions of recession and constancy of convex and concave functions are defined by Rockefellar (13, p. 69). We apply these notions to log likelihoods of full exponential families. Proofs of all theorems and corollaries are given in the appendix.

**Theorem 1.** *For some vector $\delta$ and for a full exponential family with log likelihood (3), natural parameter space $\Theta$, convex support $C$, natural statistic $Y$, and observed value of the natural statistic $y$ such that $y \in C$, the following are equivalent.*

(a) *There exists a $\theta \in \Theta$ such that $s \mapsto l(\theta + s\delta)$ is not a strictly concave function on the interval where it is finite.*
(b) *For all $\theta \in \Theta$ the function $s \mapsto l(\theta + s\delta)$ is constant on $\mathbb{R}$.*
(c) *The parameter values $\theta$ and $\theta + s\delta$ correspond to the same probability distribution for some $\theta \in \Theta$ and some $s \neq 0$.*
(d) *The parameter values $\theta$ and $\theta + s\delta$ correspond to the same probability distribution for all $\theta \in \Theta$ and all real $s$.*

(e) $\langle Y - y, \delta \rangle = 0$ *almost surely for some distribution in the family.*

(f) $\langle Y - y, \delta \rangle = 0$ *almost surely for all distributions in the family.*

(g) $\delta \in N_C(y)$ *and* $-\delta \in N_C(y)$.

(h) $\langle w, \delta \rangle = 0$, *for all* $w \in T_C(y)$.

Any vector $\delta$ that satisfies any one of the conditions of the theorem (and hence all of them) is called a *direction of constancy* of the log likelihood. The set of all directions of constancy is called the *constancy space* of the log likelihood. It is clear from (e) or (h) of the theorem that the constancy space is a vector subspace.

**Corollary 2.** *For a full exponential family, suppose $\hat{\theta}_1$ and $\hat{\theta}_2$ are maximum likelihood estimates. Then $\hat{\theta}_1 - \hat{\theta}_2$ is a direction of constancy.*

From the corollary and (d) of the theorem, we see that directions of constancy do not cause any problem for statistical inference, because all maximum likelihood estimates correspond to the same probability distribution. Thus we have uniqueness where it is important. Nonuniqueness of the MLE for the natural parameter is, at worst, merely a computational nuisance.

A family is said to be *minimal* if it has no directions of constancy. This can always be arranged by reparametrization (2, pp. 111–116; 3, pp. 13–16; see also 5, Section 1.5). The R function `glm` always uses a minimal parametrization, dropping predictors to obtain a full rank model matrix. However, insisting on minimality can complicate theoretical issues. Better to keep options open and allow for directions of constancy.

**Theorem 3.** *For some vector $\delta$ and for a full exponential family with log likelihood (3), natural parameter space $\Theta$, convex support $C$, natural statistic $Y$, and observed value of the natural statistic $y$ such that $y \in C$, the following are equivalent.*

(a) *There exists a $\theta \in \Theta$ such that the function $s \mapsto l(\theta + s\delta)$ is nondecreasing on $\mathbb{R}$.*

(b) *For all $\theta \in \Theta$ the function $s \mapsto l(\theta + s\delta)$ is nondecreasing on $\mathbb{R}$.*

(c) $\langle Y - y, \delta \rangle \leq 0$ *almost surely for some distribution in the family.*

(d) $\langle Y - y, \delta \rangle \leq 0$ *almost surely for all distributions in the family.*

(e) $\delta \in N_C(y)$.

(f) $\langle w, \delta \rangle \leq 0$, *for all* $w \in T_C(y)$.

Any vector $\delta$ that satisfies any one of the conditions of the theorem (and hence all of them) is called a *direction of recession* of the log likelihood. From now on we will simply say direction of recession or constancy to refer to directions of recession or constancy of the log likelihood. Note that every direction of constancy is a direction of recession.

**Theorem 4.** *For a full exponential family with convex support $C$ and observed value of the natural statistic $y$ such that $y \in C$, the following are equivalent.*

(a) *The MLE exists.*

(b) *Every direction of recession is a direction of constancy.*

(c) $N_C(y)$ *is a vector subspace.*

(d) $T_C(y)$ *is a vector subspace.*

This theorem provides a complete geometric solution to the problem of when the MLE exists in a full exponential family. It is the basis of the complete computational solution given in Section 3.12.

**Corollary 5.** *For a full exponential family with log likelihood* (3), *natural parameter space* $\Theta$, *convex support* $C$, *and observed value of the natural statistic* $y$ *such that* $y \in C$, *if* $\delta$ *is a direction of recession that is not a direction of constancy, then for all* $\theta \in \Theta$ *the function* $s \mapsto l(\theta + s\delta)$ *is strictly increasing on the interval where it is finite.*

### 3.4. Limits in Directions of Recession

**Theorem 6.** *For a full exponential family having log likelihood* (3), *densities* (4), *natural statistic* $Y$, *observed value of the natural statistic* $y$ *such that* $y$ *is in the convex support, and natural parameter space* $\Theta$, *if* $\delta$ *is a direction of recession,*

$$H = \{\, w \in \mathbb{R}^p : \langle w - y, \delta \rangle = 0 \,\}, \tag{9}$$

*and* $\mathrm{pr}(Y \in H) > 0$ *for some distribution in the family, and hence for all, then for all* $\theta \in \Theta$

$$\lim_{s \to \infty} f_{\theta + s\delta}(\omega) = \begin{cases} 0, & \langle Y(\omega) - y, \delta \rangle < 0 \\ f_\theta(\omega)/\mathrm{pr}_\theta(Y \in H), & \langle Y(\omega) - y, \delta \rangle = 0 \\ +\infty, & \langle Y(\omega) - y, \delta \rangle > 0 \end{cases} \tag{10}$$

*If* $\delta$ *is not a direction of constancy, then* $s \mapsto \mathrm{pr}_{\theta + s\delta}(Y \in H)$ *is continuous and strictly increasing, and* $\mathrm{pr}_{\theta + s\delta}(Y \in H) \to 1$ *as* $s \to \infty$.

We note three things about the right-hand side of (10). First, it is a probability density with respect to the distribution having parameter value $\psi$. The set where it is $+\infty$ has probability zero by Theorem 3 (d), so this is not a problem. Second, it is the density of the conditional distribution given the event $Y \in H$ of the distribution having parameter value $\theta$. Third, by Scheffé's lemma (12, p. 351) pointwise convergence of densities implies convergence in total variation, which implies convergence in distribution. Denote the right-hand side of (10) by $f_\theta(\omega \mid Y \in H)$.

It is clear that the family

$$\{\, f_\theta(\,\cdot\, \mid Y \in H) : \theta \in \Theta \,\} \tag{11}$$

is an exponential family with the same natural statistic and natural parameter as the original family. It need not be full. The natural parameter space of the full family containing it is at least as large as

$$\Theta + \Gamma_{\lim} = \{\, \theta + \gamma : \theta \in \Theta \text{ and } \gamma \in \Gamma_{\lim} \,\}, \tag{12}$$

where $\Theta$ is the natural parameter space of the original family and $\Gamma_{\lim}$ is the constancy space of the family (11). We will assume that (12) is the natural parameter space of the full family containing (11), and we will call this full family the *limiting conditional model* (LCM).

It is clear that the log likelihood for (11)

$$l_H(\theta) = \langle y, \theta \rangle - c(\theta) - \log \mathrm{pr}_\theta(Y \in H)$$

satisfies

$$l(\theta) < l_H(\theta), \qquad \theta \in \Theta.$$

Thus, if an MLE exists for the LCM, then it maximizes the likelihood in the family that is the union of the LCM and the original family. When this happens, we say we have found an MLE in the Barndorff-Nielsen completion of the original family.

### *3.5. Convex Polyhedra*

A set is a *convex polyhedron* if it is the intersection of a finite collection of closed half-spaces or, alternatively, if it is the convex hull of a finite set of points and directions (13, Chapter 19). The equivalence of these two characterizations is called the Minkowski-Weyl theorem (13, Theorem 19.1). The R function `scdd` in the contributed package `rcdd` (9) converts between these two representations of a convex polyhedron, which it calls the H-representation and V-representation, respectively.

We will need more details of V-representations. They characterize a convex polyhedron as the set of all linear combinations

$$\sum_{i \in E \cup I} b_i \alpha_i, \tag{13}$$

where the $\alpha_i$ are vectors, the $b_i$ are scalars, $E$ and $I$ are disjoint finite sets, and the $b_i$ satisfy

$$b_i \geq 0, \qquad i \in E \cup I \tag{14a}$$

and if $I$ is nonempty

$$\sum_{i \in I} b_i = 1. \tag{14b}$$

The sets $P = \{\, \alpha_i : i \in I \,\}$ and $D = \{\, \alpha_i : i \in E \,\}$ are the called points and directions, respectively, constituting the V-representation. When $P$ is empty, we write

$$C = \mathrm{con}(\mathrm{pos}\, D)$$

to denote the relationship between $D$ and the convex polyhedron $C$ it characterizes. This notation follows Rockefellar and Wets (14, Sections 2E and 3G).

When $C$ is a convex polyhedron, $N_C(y)$ and $T_C(y)$ are also convex polyhedra for each $y \in C$ and are given in terms of the H-representation of $C$ by simple formulas ([14], Theorem 6.46). Moreover, the closure operation in (7) is unnecessary when $C$ is a convex polyhedron. If a convex cone, such as $T_C(y)$ or $N_C(y)$, is polyhedral, then its V-representation can consist of directions only, so can be of the form con(pos $V$) for some finite set $V$.

For all exponential families we deal with, we assume that the convex support is polyhedral.

### *3.6. Generic Directions of Recession*

The *relative interior* of a convex set $C$, denoted rint $C$, is its interior relative to its affine hull ([13], Chapter 6). We say a vector $\delta$ is a *generic direction of recession* (GDOR) if $\delta \in$ rint $N_C(y)$ and $N_C(y)$ is not a vector subspace, where $C$ is the convex support and $y$ an observed value of the natural statistic such that $y \in C$. Since the relative interior is always nonempty ([13], Theorem 6.2), a GDOR exists if and only if none of the conditions of Theorem 4 hold.

**Theorem 7.** *For a full exponential family having polyhedral convex support $C$ and observed value of the natural statistic $y$ such that $y \in C$, let $T_C(y) =$ con(pos $V$), and define*

$$L = \{ v \in V : -v \in T_C(y) \}.$$

*Then a generic direction of recession exists if and only if $L \neq V$, in which case a vector $\delta$ is a generic direction of recession if and only if*

$$\langle w, \delta \rangle = 0, \qquad w \in L \tag{15a}$$

$$\langle w, \delta \rangle < 0, \qquad w \in V \setminus L \tag{15b}$$

**Corollary 8.** *Under the assumptions of the theorem, a generic direction of recession is not a direction of constancy.*

If $B$ is a set of vectors, let span $B$ denote the smallest vector subspace containing $B$. Also for any vector $x$, let $x +$ span $B = \{ x + v : v \in$ span $B \}$.

**Corollary 9.** *Under the assumptions of the theorem, suppose $\delta$ is a generic direction of recession, and $H$ is defined by (9). Then $T_{C \cap H}(y) =$ span $L$, and $C \cap H = C \cap (y + $ span $L)$.*

The theorem and corollaries explain the purpose of generic directions of recession. By Corollary 8, a GDOR implies the MLE does not exist in the conventional sense, so we seek it in the LCM described in Section 3.4 using the GDOR as the $\delta$ in Theorem 6. Suppose that $C \cap H$ is the convex support of the LCM. Then $T_{C \cap H}(y)$ being a vector subspace implies that the MLE in the LCM exists by Theorem 4 (c). Thus finding one GDOR allows us to find the MLE in the Barndorff-Nielsen completion.

### 3.7. Assumptions

We summarize the assumptions we have made above. We deal with a full exponential family having convex support $C$ that is a polyhedron. If every direction of recession is a direction of constancy, then we need no further assumptions. Otherwise, let $\delta$ be a generic direction of recession, let $C$ be the convex support, let $Y$ be the natural statistic, let $y$ be an observed value of the natural statistic satisfying $y \in C$, and let $H$ be defined by (9). We assume the event $Y \in H$ has positive probability so the LCM defined in Section 3.4 exists. We further assume that $C \cap H$ is the convex support of this LCM, so that by Corollary 9 the MLE in this LCM exists. We further assume that the natural parameter space of this LCM is given by (12), so that confidence intervals (Section 3.16 below) work.

### 3.8. Comparison with Pre-Existing Theory

The pre-existing theory of the Barndorff-Nielsen completion (2; 3; 5) says the MLE lies in the LCM whose convex support, what we are calling $C \cap H$, is the unique face of $C$ containing $y$ in its relative interior (5, Chapter 4 generalizes this). Thus the pre-existing approach makes it clear that the LCM containing the MLE is unique and does not depend on the GDOR, which is in general not unique. In our approach, uniqueness comes from the assertion $C \cap H = C \cap (y + \text{span } L)$ in Corollary 9. This makes it clear that, although the hyperplane $H$ does depend on the GDOR $\delta$ used to define it, the convex support $C \cap H$ of the LCM does not depend on $H$, hence does not depend on $\delta$.

### 3.9. Natural Affine Submodels

In most applications of exponential family theory, we start with a very large exponential family, which we call *saturated* and which has too many parameters to estimate well. Then we consider *natural affine submodels*, parametrized by

$$\theta = a + M\beta,$$

where $\theta$ is the natural parameter of the saturated model, $\beta$ is the natural parameter of the natural affine submodel, $a$ is a known vector, and $M$ is a known matrix. In the terminology of the R function `glm`, $a$ is called the *offset vector* and $M$ is called the *model matrix*.

Because

$$\langle y, a + M\beta \rangle = \langle y, a \rangle + \langle M^T y, \beta \rangle,$$

where the two bilinear forms on the right-hand side have different dimensions, and the first term on the right-hand side does not contain the parameter and can be dropped from the log likelihood, the submodel is itself an exponential family with natural statistic $M^T y$ and natural parameter $\beta$. Thus everything said above applies to natural affine submodels, we just work with the convex support of $M^T Y$ rather than of $Y$.

### *3.10. Tangent Cones of Models and Affine Submodels*

Let $C_{\mathrm{sat}}$ denote the convex support of the saturated model and $C_{\mathrm{sub}}$ that of the natural affine submodel. By Theorems 6.43 and 6.46 in (14),

$$T_{C_{\mathrm{sub}}}(M^T y) = \mathrm{cl}\{\, M^T w : w \in T_{C_{\mathrm{sat}}}(y) \,\} \tag{16}$$

and the closure operation is not necessary if $C_{\mathrm{sat}}$ is polyhedral. Moreover, it is clear that if $T_{C_{\mathrm{sat}}}(y) = \mathrm{con}(\mathrm{pos}\, V_{\mathrm{sat}})$, then $T_{C_{\mathrm{sub}}}(y) = \mathrm{con}(\mathrm{pos}\, V_{\mathrm{sub}})$, where

$$V_{\mathrm{sub}} = \{\, M^T w : w \in V_{\mathrm{sat}} \,\}. \tag{17}$$

### *3.11. Tangent Cones of Saturated Models*

If we assume nothing about $C_{\mathrm{sat}}$ except that it is polyhedral, then computation of a $V$-representation $V_{\mathrm{sat}}$ for it can be arbitrarily difficult. Thus we shall say nothing about the general case and focus on cases of practical interest. For generalized linear models, the convex support of the saturated model is a Cartesian product. For logistic regression each component of the response vector is Bernoulli and $C_{\mathrm{sat}} = [0, 1]^p$. For Poisson regression, each component of the response vector is Poisson and $C_{\mathrm{sat}} = [0, \infty)^p$.

When $C_{\mathrm{sat}}$ is a Cartesian product, $T_{C_{\mathrm{sat}}}(y)$ can be calculated coordinatewise (14, Proposition 6.41). Let $e_i$ denote the unit vector in the $i$-th coordinate direction (every coordinate is zero except for the $i$-th, which is one). Then $e_i$ is a tangent vector at $y$ if $y_i$ is not at the upper bound of its range, and $-e_i$ is a tangent vector at $y$ if $y_i$ is not at the lower bound. Hence when $C_{\mathrm{sat}}$ is a Cartesian product, its V-representation $V_{\mathrm{sat}}$ contains $e_i$ or $-e_i$ or both for each $i$.

In log-linear models for categorical data analysis, if Poisson response is assumed, then $C_{\mathrm{sat}} = [0, \infty)^p$, just as in Poisson regression. If multinomial or product-multinomial response is assumed, then $C_{\mathrm{sat}}$ is not a Cartesian product, but, as we shall see (Section 3.17 below), the MLE and GDOR are the same regardless of the assumed response distribution. Hence the solution for the Poisson response case can be used for the other cases.

Thus this analysis of the Cartesian product case suffices for the vast majority of applications, and it suffices for our examples. An application where the convex support of the saturated model is not a Cartesian product and this analysis would not suffice is unconditional aster models (11).

### *3.12. Calculating the Linearity*

Next we determine the *linearity* of $V_{\mathrm{sub}}$

$$L_{\mathrm{sub}} = \{\, w \in V_{\mathrm{sub}} : -w \in \mathrm{con}(\mathrm{pos}\, V_{\mathrm{sub}}) \,\}. \tag{18}$$

This sounds like a complicated operation, and it is, but the `rcdd` package has a function `linearity` that does it by repeated linear programming.

Having found the linearity, we have a complete computational solution to the problem of when the MLE exists in a full exponential family with polyhedral convex support. It exists if and only if $L_{\mathrm{sub}} = V_{\mathrm{sub}}$.

All functions in the `rcdd` package use two forms of arithmetic. One is the default computer arithmetic used by all other R functions. Answers produced using that arithmetic are inexact, so one is uncertain whether the $L_{\mathrm{sub}}$ produced is actually correct. The other form of arithmetic is exact, infinite-precision, rational arithmetic. Answers produced using that arithmetic are exact, so one is certain that the $L_{\mathrm{sub}}$ produced is actually correct, but only if the vectors in $V_{\mathrm{sub}}$ are also produced exactly using either integer arithmetic or rational arithmetic.

According to comments in the source code (starting at line 3062 of the file **cddlp.c** of the source code for the `rcdd` package (9), which comes from the **cddlib** library (4), for each $w \in V_{\mathrm{sub}}$ the R function `linearity` solves the linear programming problem

$$
\begin{aligned}
&\text{maximize} \\
&\langle w, \delta \rangle \\
&\text{subject to} \\
&\langle v, \delta \rangle \geq 0, \qquad v \in V_{\mathrm{sub}} \setminus \{w\}
\end{aligned}
\tag{19}
$$

where $\delta$ is the state vector of the linear programming problem.

**Theorem 10.** *A vector $w$ is in the linearity* (18) *if and only if the optimal value of the linear program* (19) *is nonpositive.*

### 3.13. Calculating Generic Directions of Recession

If $L_{\mathrm{sub}} \neq V_{\mathrm{sub}}$, then $T_{C_{\mathrm{sub}}}(M^T y)$ is not a vector subspace, hence there exists a generic direction of recession. By Theorem 7, $\delta$ is a GDOR if and only if

$$
\begin{aligned}
\langle w, \delta \rangle = 0, \qquad w \in L_{\mathrm{sub}} && \text{(20a)} \\
\langle w, \delta \rangle < 0, \qquad w \in V_{\mathrm{sub}} \setminus L_{\mathrm{sub}} && \text{(20b)}
\end{aligned}
$$

Hence we can find one such $\delta$ by solving the following linear program

$$
\begin{aligned}
&\text{maximize} \\
&\epsilon \\
&\text{subject to} \\
&\epsilon \leq 1 \\
&\langle v, \delta \rangle = 0, \qquad v \in L_{\mathrm{sub}} \\
&\langle v, \delta \rangle \leq -\epsilon, \qquad v \in V_{\mathrm{sub}} \setminus L_{\mathrm{sub}}
\end{aligned}
$$

where $\delta$ is a $q$-dimensional vector, $\epsilon$ is a scalar, and $(\delta, \epsilon)$ is the state vector of the linear program (so the dimension is $q + 1$). The $\delta$ part of the solution is a generic direction of recession. The $\epsilon$ part does not matter.

The idea for using this particular linear program came from the documentation for the `dd_ExistsRestrictedFace2` function in the `cddlib` library (written by K. Fukuda). The `rcdd` package does not provide an interface to this `cddlib` function, but it does provide a function `lpcdd` that does linear programming and can be used to solve this linear program.

### *3.14. Calculating Maximum Likelihood Estimates*

#### *3.14.1. In the Original Family*

There is little to be said about calculating the MLE in the original family. When we have found that $L_{\text{sub}} = V_{\text{sub}}$, then we know the MLE exists and can use available software to find it. We will use the R function `glm` for our examples.

There is one issue worth mentioning. If the model is non-identifiable, so the MLE is non-unique, the R function `glm` is smart enough to drop enough predictors to produce an identifiable model. However, its method of doing so is not guaranteed because of inexactness of the default computer arithmetic.

We can use the `redundant` function in the `rcdd` package applied to the columns of the model matrix $M$ to reduce to a linearly independent subset. If $M$ was calculated using integer or rational arithmetic, and `redundant` uses rational arithmetic, then this operation will be exact.

#### *3.14.2. In the Completion*

When we have found that $L_{\text{sub}} \neq V_{\text{sub}}$ and have found a generic direction of recession $\delta$, we still need to characterize the support of the LCM. We can characterize this two ways using either of

$$H_{\text{sub}} = \{\, w \in \mathbb{R}^q : \langle w - M^T y, \delta \rangle = 0 \,\}$$
$$H_{\text{sat}} = \{\, w \in \mathbb{R}^p : \langle w - y, M\delta \rangle = 0 \,\}$$

where $p$ and $q$ are the dimensions of the saturated model and affine submodel, respectively. Then the LCM conditions on the event $M^T Y \in H_{\text{sub}}$ or $Y \in H_{\text{sat}}$, which is the same event characterized two different ways, the latter usually simpler.

**Theorem 11.** *In the setup of Sections 3.12 through 3.14, define*

$$L_{\text{sat}} = \{\, v \in V_{\text{sat}} : M^T v \in L_{\text{sub}} \,\}.$$

*Then the support of the limiting conditional model is*

$$C_{\text{sat}} \cap H_{\text{sat}} = C_{\text{sat}} \cap (y + \text{span}\, L_{\text{sat}}). \tag{21}$$

In the case where the convex support of the saturated model is a Cartesian product, the support of the LCM simply constrains $Y_i = y_i$ for $i$ such that $e_i \notin L_{\text{sat}}$, that is, the $i$-th component of the response is unconstrained if $e_i \in L_{\text{sat}}$ and is constrained to be equal to its observed value if $e_i \notin L_{\text{sat}}$.

Our analysis of maximum likelihood estimation in the Barndorff-Nielsen completion is now finished. At least in the Cartesian product case, the maximum likelihood problem in the completion is of the same form as the original problem. The only difference is that we constrain certain components of the response vector to their observed values. This can be achieved by removing those components from the response vector and proceeding as if the resulting subvector were the entire response vector. If, for example, we are using the R function `glm` to fit models, we merely delete certain elements of the response vector and the corresponding rows of the model matrix (or the data frame containing the data if we are using a formula to specify the model) and proceed normally. This LCM produced by deleting some components of the response will always be non-identifiable, because $\delta$ will always be a direction of constancy for the LCM and there may be other directions of constancy. The `glm` function, however, can deal with this issue. Furthermore, even in the rare case when the `glm` function may be confused, we can find a full rank model matrix having the same column space as the original model matrix using the function `redundant` in the `rcdd` package, as described in the preceding section.

### 3.15. Likelihood Ratio Tests

Given two nested natural affine submodels, the maximum value of the log likelihood can be calculated for each submodel by available software, such as the R function `glm`, which goes uphill on the log likelihood until reaching a point where the log likelihood is nearly flat, in which case the value of the log likelihood is nearly the maximum. If the MLE does not exist in the conventional sense, then the natural parameter estimates will be large but not infinite, and the `glm` function may or may not give a warning about lack of convergence. If the MLE does not exist in the conventional sense, then the natural parameter estimates are infinitely wrong, but the value of the maximized log likelihood is nearly correct. Thus we can correctly calculate the likelihood ratio test statistic without using any of the theory developed here.

This does no good, however, because the usual asymptotics of the likelihood ratio test (Wilks' theorem) do not hold in the case where the MLE for the null model does not exist in the conventional sense. In this case, the following simple correction, suggested by S. Fienberg (personal communication) seems reasonable. Apply the theory developed here to the null model, determining $L_{\text{sat}} \neq V_{\text{sat}}$. Let $M_0$ and $M_1$ be the model matrices for the null and alternative natural affine submodels.

If we apply Wilks' theorem to the LCM for the null hypothesis, we obtain the result that the deviance (twice the log likelihood ratio) is approximately chi-squared with degrees of freedom which is the difference in dimension of $M_1^T(\text{span } L_{\text{sat}})$ and of $M_0^T(\text{span } L_{\text{sat}})$. Assuming that $M_0$, $M_1$, and $L_{\text{sat}}$ were determined exactly using rational arithmetic, the degrees of freedom can be determined exactly by applying the `redundant` function in the `rcdd` package to the sets $\{\, M_i^T w : w \in L_{\text{sat}} \,\}$, $i = 0, 1$.

This asymptotic approximation may or may not hold depending on the sample size and on how close the observed value of the natural statistic is to the boundary of the convex support of the LCM. By construction, it cannot be on the boundary, but if it is close the asymptotic approximation can be bad. If one is worried about the validity of the asymptotic approximation, one can always do a parametric bootstrap calculation based on the LCM for the null hypothesis.

### *3.16. Confidence Intervals*

Confidence intervals are more complicated than hypothesis tests. Confidence intervals for both natural and mean value parameters are of interest. The R function `predict.glm` provides either, depending on the value of its `type` argument. We will also provide either.

Before getting into details, we should first note that confidence intervals are often inappropriate from a theoretical point of view. When there is not a single scalar parameter of interest, a confidence region for the vector parameter of interest should, theoretically, be provided. However, high-dimensional confidence regions are unvisualizable and of no interest to users. Thus statisticians usually provide something users think they can interpret, which is multiple confidence intervals, often not adjusted for simultaneous coverage. That is what, for example, `predict.glm` provides. We will follow the usual practice, providing such confidence intervals in our examples. The accompanying technical report (7) discusses confidence regions.

In the case where no GDOR exists and the MLE exists in the conventional sense, we say the original family is the LCM, which it is when we take $\delta = 0$ in Theorem 6. Then there is always an LCM. We calculate confidence intervals conditional on the LCM. If we achieve approximately the desired confidence level conditionally, then we also achieve it unconditionally.

### *3.16.1. In the Limiting Conditional Model*

Our strategy is, like it was with hypothesis tests, to use standard methods applied to the LCM as much as possible. Unlike it was with hypothesis tests, however, no simple modification of standard methods tells us all we want to know. When we fit the LCM, available software gives us confidence intervals for its regression coefficients, for components of its linear predictor vector, and for its mean value parameters. Its mean value parameter vector is restricted to its convex support $C \cap H$. No procedure based on the LCM can indicate how close the unknown true distribution of the data is to being concentrated on $C \cap H$. To accomplish that, we need to refer to the original model.

### *3.16.2. In the Original Model*

Let $\hat{\theta}$ be an MLE in the LCM and let $\Gamma_{\text{lim}}$ be the constancy space of the LCM, so that $\hat{\theta} + \gamma$ is also an MLE in the LCM for $\gamma \in \Gamma_{\text{lim}}$. Let $\delta$ be a GDOR. We know from Theorem 6 that, regardless of the $\gamma$ chosen, distributions in the original model having parameter values $\hat{\theta} + \gamma + s\delta$ converge to the MLE distribution in the LCM as $s \to \infty$. We seek an answer to the question how close $s$ must be to infinity. Thus we seek one-sided intervals for the scalar parameter $s$ of the form $(\hat{s}_\gamma, \infty)$.

The conventional conservative $P$-value for the upper-tailed test having test statistic $\langle Y, \delta \rangle$, null hypothesis $\hat{\theta} + \gamma + s\delta$, and observed data $y$ is

$$\text{pr}_{\hat{\theta} + \gamma + s\delta}(\langle Y - y, \delta \rangle \geq 0). \tag{22}$$

A conventional level $\alpha$ test rejects when (22) is less than or equal to $\alpha$. A conservative $1 - \alpha$ confidence interval consists of the set of $s$ values that are not rejected at level $\alpha$. In the only case of interest to us, where $\langle y, \delta \rangle$ is the largest possible value of $\langle Y, \delta \rangle$ so the event $\langle Y - y, \delta \rangle \geq 0$ is the same as $Y \in H$ up to a set of measure zero,

$$\hat{s}_\gamma = \inf\{\, s \in \mathbb{R} : \hat{\theta} + \gamma + s\delta \in \Theta \text{ and } \text{pr}_{\hat{\theta} + \gamma + s\delta}(Y \in H) > \alpha \,\}.$$

Since $s \mapsto \text{pr}_{\hat{\theta} + \gamma + s\delta}(Y \in H)$ is continuous and strictly increasing by Theorem 6, usually $\hat{s}_\gamma$ is the unique $s$ such that $\text{pr}_{\hat{\theta} + \gamma + s\delta}(Y \in H) = \alpha$. The accompanying technical report (7) explains how these intervals can be made exact using the method of fuzzy confidence intervals (8).

Since all of the intervals $(\hat{s}_\gamma, \infty)$ for different $\gamma$ are based on the same test statistic $\langle Y, \delta \rangle$, they have simultaneous coverage at least $1 - \alpha$. Hence the union

$$\{\, \hat{\theta} + \gamma + s\delta : \gamma \in \Gamma_{\text{lim}} \text{ and } s > \hat{s}_\gamma \,\} \tag{23}$$

of those natural parameter values provides a confidence region for natural parameter values showing "how close to infinity" they may be.

This completes what the original model can say about natural parameters. When converted to statements about mean value parameters, these confidence intervals say how close the true unknown distribution may be to being concentrated on $C \cap H$.

### *3.16.3. The Cartesian Product Case*

In the case where the convex support of the saturated model is a Cartesian product, all of this simplifies somewhat. We know that a confidence region for the mean value parameters of the LCM only involves response variables that are not constrained to be at their observed values in the LCM. We take confidence intervals, such as those provided by the R function `predict.glm` applied to the limiting conditional model, to be adequate for describing those components of the response.

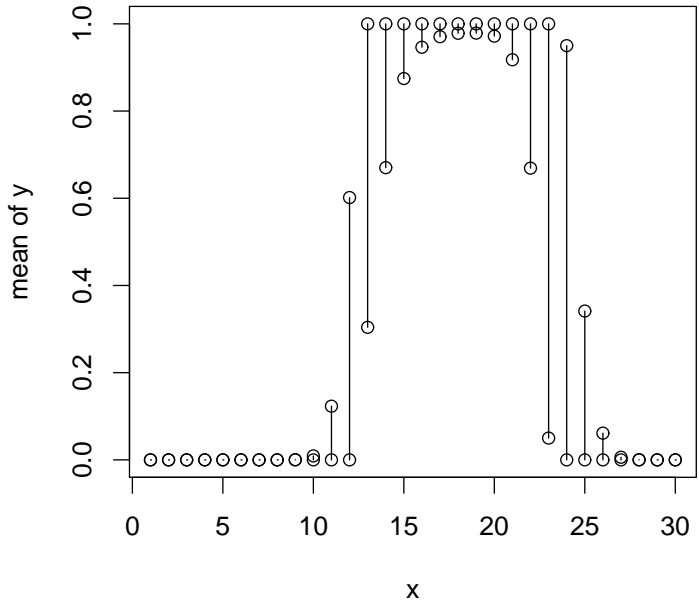

FIG 2. *Careful confidence intervals for the mean value parameter in the example of Section 2.2. Outer points having values zero or one are MLE for the mean value parameter, which are also the components of the observed data vector. Inner points are mean values corresponding to $\hat{\beta} + \gamma + \hat{s}_\gamma \delta$, where $\hat{s}_\gamma$ is the lower bound for a 95% one-sided confidence interval for s, where $\hat{\beta} = 0$ and $\delta$ is given by (1), and where $\gamma$ is chosen differently for each mean value to make the interval as large as possible. Compare Figure 1.*

Our one-sided intervals come into play in computing one-sided confidence intervals for the mean value parameters of the other components of the response. Since the MLE of their mean value parameters are on the boundary, one-sided intervals are the only kind that make sense. We distinguish two cases.

The first case is where $\Gamma_{\text{lim}} = \{\, s\delta : s \in \mathbb{R} \,\}$. In this case, all intervals $(\hat{s}_\gamma, \infty)$ for different $\gamma \in \Gamma_{\text{lim}}$ correspond to the same natural parameter values. Hence one such interval adequately describes the variability in mean value parameters for components of the response having MLE at the boundary. This case contains our contingency table example (Section 2.3).

The second case is where $\Gamma_{\text{lim}}$ contains vectors other than scalar multiples of the GDOR. In this case, intervals $(\hat{s}_\gamma, \infty)$ for different $\gamma \in \Gamma_{\text{lim}}$ correspond to different natural parameter values. Since they have simultaneous coverage, we can use the union (23) to determine these confidence intervals. Figure 2 shows confidence intervals for mean value parameters based on this idea for our logistic regression example (Section 2.2).

Calculation of these intervals theoretically requires knowing $\hat{s}_\gamma$ for all $\gamma \in \Gamma_{\text{lim}}$, which entails an infinite amount of work. We have only calculated for a large finite set of $\gamma$.

### 3.16.4. Calculating the Constancy Space

By Theorems 1 and 4 the constancy space is $N_C(y)$ in the case where every direction of recession is a direction of constancy. By Corollary 9 the tangent space $T_{C \cap H}(y)$ in the LCM is span $L$. Hence the constancy space is

$$N_{C \cap H}(y) = \{ \delta \in \mathbb{R}^p : \langle v, \delta \rangle = 0, \ v \in \text{span}\, L \}.$$

In the case of a natural affine submodel this becomes

$$\Gamma_{\text{lim}} = N_{C_{\text{sub}} \cap H_{\text{sub}}}(M^T y) = \{ \delta \in \mathbb{R}^q : \langle v, \delta \rangle = 0, \ v \in \text{span}\, L_{\text{sub}} \}. \tag{24}$$

Since $L_{\text{sub}}$ is a V-representation of its span, a call to the function `scdd` in the `rcdd` package will compute an H-representation of its span which is also a V-representation for (24), that is, a basis for the constancy space of the LCM.

When calculating (24) for the purpose of generating one-sided confidence intervals it is clear that we do not need to let $\gamma$ range over the whole constancy space (24) because points $\gamma$ and $\gamma + s\delta$ lead to the same one-sided confidence intervals. Hence it is enough to use the subspace of $\Gamma_{\text{lim}}$ orthogonal to $\delta$, which is calculated by feeding $L_{\text{sub}} \cup \{\delta\}$ to the R function `scdd` for conversion to an H-representation, which will also be a basis of the desired subspace.

### 3.17. Multinomial Sampling

Consider one contingency table and one vector of observed data, but two models: Poisson sampling and multinomial sampling. As is well known (1, Section 8.6.7), the maximum likelihood estimates for the mean value parameters are the same for both sampling schemes. But much more is the same.

Suppose we consider the natural statistic to be the vector of cell counts for both models, so both have the same natural statistic and natural parameter. For Poisson sampling, there are no directions of constancy and the MLE for the natural parameter is unique. For multinomial sampling, the vector $\gamma = (1, 1, \ldots, 1)$ is a direction of constancy, and the MLE for the natural parameter is not unique. However, it is easy to see that the unique MLE for the Poisson model is also an MLE for the multinomial model.

Moreover, when we use the same natural statistic and parameter for both models, the computational geometry is similar. If subscripts $P$ and $M$ refer to the Poisson and multinomial models, respectively, then

$$C_{\text{sat},M} = \{ y \in C_{\text{sat},P} : \langle y, \gamma \rangle = n \}$$

where $n$ is the sample size. Also note that $\gamma$ is the first column of $M$, the "intercept" column. From this it follows that

$$T_{C_{\mathrm{sub},M}}(M^T y) = \{\, v \in T_{C_{\mathrm{sub},P}}(M^T y) : \langle v, e_1 \rangle = 0 \,\}$$

where $e_1 = (1, 0, \dots, 0)$. And from this it follows by Theorem 6.42 in ([14]) that

$$N_{C_{\mathrm{sub},M}}(M^T y) \supset \{\, v + s e_1 : v \in N_{C_{\mathrm{sub},P}}(M^T y), \ s \in \mathbb{R} \,\}.$$

Hence every direction of recession in the Poisson model is also one in the multinomial model. Similarly, every GDOR in the Poisson model is also one in the multinomial model. Hence the support of the LCM calculated using such a GDOR is also correct. Hence parameter estimates for the LCM are the same for both Poisson and multinomial sampling, as are the asymptotic likelihood ratio tests described in Section 3.15 above.

The only thing we need to change for multinomial sampling is our one-sided confidence intervals described in Section 3.16.2 above, because they are based on exact probabilities that differ between the two models. The modification is obvious: in calculating $\mathrm{pr}_{\hat{\theta} + \gamma + s\delta}(Y \in H)$ we use the multinomial distribution rather than the Poisson distribution. Table 2 in the accompanying technical report ([7]) shows the analog of our Table 2 modified using multinomial rather than Poisson sampling.

Product-multinomial sampling is very similar to the situation for multinomial sampling. We only note that each sum of components of the response that is fixed must be a column of the model matrix if the MLE for Poisson sampling are to match that for product-multinomial ([1], Section 8.6.7). Details are left as an exercise for the reader.

## 4. Discussion

Part of the impetus for writing this article was having to say to a scientist, "you are just out of luck, the solution is 'at infinity' and this problem is well understood theoretically but software just doesn't handle it — the only thing you can do with existing software is fit a smaller model that doesn't have the solution 'at infinity' even though this smaller model admittedly (1) does not fit the data and (2) does not address the questions of scientific interest." I am glad I will never have to say this again.

Questions of usability and user interface remain. The R function `glm` and generalized linear models software in other statistical computing environments (SCE) should just do the right thing when the MLE does not exist in the conventional sense. This would not only require additional programming — not much for R but much more for other SCE that do not have a computational geometry package — but also would break backward compatibility. The R function `predict.glm`, like other methods of the generic function `predict`, specifies confidence intervals by estimates and standard errors (components `fit` and `se.fit` of the returned object). This clearly will not do for one-sided intervals.

The theorems and corollaries of this article are very general, applying to all known applications. The theory of hypothesis tests presented (due to S. Fienberg) is also general. For confidence intervals and regions, however, open research questions remain.

Aster models (11) allow components of the response vector to be dependent. For example, the conditional distribution of $Y_j$ given $Y_k$ can be Binomial$(Y_k, p_j)$. When one has $y_j = y_k$ in the observed data, it can happen that the mean value parameter vector with components $\mu_j = E(Y_j)$ has $\hat{\mu}_j = \hat{\mu}_k$ at the MLE, in which case a one-sided confidence interval for $\mu_k - \mu_j$ is appropriate, but one-sided intervals for either $\mu_j$ or $\mu_k$ do not make sense. More generally, the R function `predict.aster`, unlike other methods of the R generic function `predict`, computes confidence intervals for arbitrary linear functions of the mean value parameter vector, which is often useful in applications (the example in 11 and all three examples in 15). It is not clear how the one-sided and two-sided confidence intervals discussed in Sections 3.16.2 and 3.16.1 above combine in this setting.

Even though issues remain, ordinary R users should be able to follow the examples in the accompanying technical report (7) to make valid hypothesis tests and confidence intervals for GLM and loglinear models for contingency tables in cases where the MLE does not exist in the conventional sense and previously available software was useless.

## Acknowledgments

Much of the theory of this article comes from the thesis (5) done under the supervision of Elizabeth A. Thompson.

## Appendix A: Proofs

*Proof of Theorem 1.* Clearly, $s \mapsto l(\theta + s\delta)$ fails to be strictly concave if and only if $s \mapsto c(\theta + s\delta)$ fails to be strictly convex, and by Theorem 2.1 in (5) this happens if and only if $\langle Y, \delta \rangle$ is concentrated at one point, in which case this point must be $\langle y, \delta \rangle$ so (e) holds. Since all distributions in the family are mutually absolutely continuous by (4), (e) implies (f), which trivially implies (e). If (f) holds, then by (5)

$$
\begin{aligned}
c(\theta + s\delta) &= c(\psi) + \log E_\psi \left( e^{\langle Y, \theta + s\delta - \psi \rangle} \right) \\
&= c(\psi) + s\langle y, \delta \rangle + \log E_\psi \left( e^{\langle Y, \theta - \psi \rangle} \right) \\
&= c(\theta) + s\langle y, \delta \rangle
\end{aligned}
\tag{25}
$$

Hence (b) holds, and (b) clearly implies (a). We have now proved that (a), (b), (e), and (f) are equivalent.

Also (25) implies (d) by (4), so (f) implies (d). Trivially, (d) implies (c). Conversely, if (c) holds, then $f_\theta$ and $f_{\theta + s\delta}$ must be equal almost surely, hence by (4)

$$
\log f_{\theta + s\delta}(\omega) - \log f_\theta(\omega) = s\langle Y(\omega), \delta \rangle - c(\theta + s\delta) + c(\theta)
$$

almost surely, hence $\langle Y, \delta \rangle$ is constant almost surely, and the constant must be $\langle y, \delta \rangle$; hence (e) holds. We have now proved that (a) through (f) are equivalent.

By definition of normal cone and convex support, (e) and (g) are equivalent, and (g) and (h) are equivalent by the polarity relationship of normal and tangent cones (14, Theorem 6.9 and Corollary 6.30). □

*Proof of Corollary 2.* By Theorem 7.1 and p. 140 in (2), $l$ is concave; thus we must have

$$l\big(t\hat{\theta}_1 - (1-t)\hat{\theta}_2\big) \geq tl(\hat{\theta}_1) + (1-t)l(\hat{\theta}_2), \qquad 0 < t < 1, \tag{26}$$

and since $\hat{\theta}_1$ and $\hat{\theta}_2$ are MLE, (26) must actually hold with equality. Thus by (a) of the theorem $\hat{\theta}_1 - \hat{\theta}_2$ is a direction of constancy. □

For the proof of Theorem 3 we use Corollary 2.4.1 in (5), which relies on Theorem 2.3 in (5), but the proof of that theorem given in (5) is murky at best. So we give a corrected version.

*Corrected Proof of Theorem 2.3 in (5).* Equation (2.5) in (5) contains an obvious typographical error. It should read

$$(\text{rc}\log c)(\phi) = \lim_{s \to \infty} \frac{\log c(\theta + s\phi) - \log c(\theta)}{s}$$
$$= \lim_{s \to \infty} \log\left(\left[\frac{c(\theta + s\phi)}{c(\theta)e^{s\sigma_K(\phi)}}\right]^{1/s} e^{\sigma_K(\phi)}\right)$$

The rest of the proof of the $\lambda(H_\phi) > 0$ case is correct. In the proof of the of the $\lambda(H_\phi) = 0$ case, the last displayed formula of the proof is incorrect. Clearly

$$e^{a - \sigma_K(\phi)} F_\theta(A)^{1/s} \to e^{a - \sigma_K(\phi)}, \qquad \text{as } s \to \infty.$$

However, since $a < \sigma_K(\phi)$ was arbitrary, the limit can be made arbitrarily close to 1, and we see that

$$\left[\frac{c(\theta + s\phi)}{c(\theta)e^{s\sigma_K(\phi)}}\right]^{1/s} \to 1, \qquad \text{as } s \to \infty,$$

as is required for the completion of the proof. □

*Proof of Theorem 3.* The equivalence of (a) and (b) is Theorem 8.6 in (13). The equivalence of (a) and (c) is Corollary 2.4.1 in (5). The equivalence of (c) and (d) is mutual absolute continuity of the distributions in an exponential family. The equivalence of (c) and (e) is immediate from our definition (8) of the normal cone. The equivalence of (e) and (f) is the polarity relationship of tangent and normal cones (14, Theorem 6.9 and Corollary 6.30). □

*Proof of Theorem 4.* That (a) and (b) are equivalent is Theorem 2.5 in (5). That (b) and (c) are equivalent follows from (g) of Theorem 1 and (e) of Theorem 3. That (c) and (d) are equivalent is the polarity relationship of tangent and normal cones. □

*Proof of Corollary 5.* By assumption $s \mapsto l(\theta + s\delta)$ is a nondecreasing function. Suppose to get a contradiction that

$$l(\theta + s_1\delta) = l(\theta + s_2\delta) \tag{27}$$

for some $s_1$ and $s_2$ such that both sides of (27) are finite and $s_1 < s_2$. In order that $l$ be nondecreasing we must have

$$l(\theta + s_1\delta) = l(\theta + s\delta), \qquad s_1 \leq s \leq s_2$$

but then $\delta$ is a direction of constancy by Theorem 1 (a). □

*Proof of Theorem 6.* Except for the last sentence, this follows immediately from Theorem 2.2 in (5). From (4)

$$
\begin{aligned}
\mathrm{pr}_{\theta+s\delta}(Y \in H) &= E_\psi\big\{ I_H e^{\langle Y, \theta+s\delta-\psi\rangle - c(\theta+s\delta) + c(\psi)} \big\} \\
&= e^{s\langle y,\delta\rangle - c(\theta+s\delta) + c(\theta)} E_\psi\big\{ I_H e^{\langle Y, \theta-\psi\rangle - c(\theta) + c(\psi)} \big\} \\
&= e^{s\langle y,\delta\rangle - c(\theta+s\delta) + c(\theta)} \, \mathrm{pr}_\theta(Y \in H)
\end{aligned}
$$

where $I_H$ denotes the indicator function of the event $Y \in H$. By Corollary 5, the function $s \mapsto \langle y, \theta + s\delta \rangle - c(\theta + s\delta)$ is strictly increasing, hence so is $s \mapsto \mathrm{pr}_{\theta+s\delta}(Y \in H)$. That $\mathrm{pr}_{\theta+s\delta}(Y \in H) \to 1$ as $s \to \infty$ follows from Scheffé's lemma (see the comments following the theorem). The continuity assertion follows from the fact that the moment generating function of the random variable $\langle Y, \delta \rangle$ is

$$
\begin{aligned}
E_\theta\big\{ e^{s\langle Y,\delta\rangle} \big\} &= E_\psi\big\{ e^{\langle Y, \theta+s\delta-\psi\rangle - c(\theta) + c(\psi)} \big\} \\
&= e^{c(\theta+s\delta) - c(\theta)}
\end{aligned}
$$

Hence $s \mapsto c(\theta + s\delta)$ is actually infinitely differentiable and so is $s \mapsto \mathrm{pr}_{\theta+s\delta}(Y \in H)$. □

*Proof of Theorem 7.* Suppose $L = V$. Then $\mathrm{con}(\mathrm{pos}\, V)$ is the subspace spanned by $V$, in which case a GDOR does not exist by Theorem 4.

Suppose $L \neq V$. Then by the polarity relationship of normal and tangent cones for each $v \in V \setminus L$ there exists $\delta_v \in N_C(y)$ such that $\langle v, \delta_v \rangle < 0$. Hence $-\delta_v \notin N_C(y)$ and $N_C(y)$ is not a vector subspace. So a GDOR does exist by Theorem 4.

Let $\delta^* = \sum_{v \in V \setminus L} \delta_v$. Then $\delta^*$ satisfies (15a) and (15b). Observe that $\delta \in N_C(y)$ if and only if (15a) holds and (15b) holds with $<$ replaced by $\leq$. Then it is clear that for every $\delta \in N_C(y)$ there exists $t > 1$ such that $t\delta^* + (1-t)\delta$ is in $N_C(y)$. Hence $\delta^* \in \mathrm{rint}\, N_C(y)$ by Theorem 6.4 in (13). It now follows from Proposition 2.42 in (14) that the set of points satisfying (15a) and (15b) is $\mathrm{rint}\, N_C(y)$. □

*Proof of Corollary 8.* In the proof of the theorem we saw that if a GDOR exists, then $L \neq V$ and $N_C(y)$ is not a vector subspace. □

*Proof of Corollary 9.* Since $C$ is polyhedral convex, every tangent vector is of the form $s(w - y)$ for some $w \in C$ and $s \geq 0$, that is, the closure operation in (7) is not necessary. This implies, in particular, that for each $v \in L$ there exist points $w_{v,+}$ and $w_{v,-}$ in $C$ and positive scalars $s_{v,+}$ and $s_{v,-}$ such that $\pm v = s_{v,\pm}(w_{v,\pm} - y)$. Observe that these $w_{v,\pm}$ are also in $C \cap H$, but no $w \in V \setminus L$ is in $C \cap H$. Thus $T_{C \cap H}(y) = \mathrm{con}(\mathrm{pos}\, L) = \mathrm{span}\, L$. Since $y + \mathrm{span}\, L \subset H$, we have $C \cap H \supset C \cap (y + \mathrm{span}\, L)$. If $C \cap H \not\subset C \cap (y + \mathrm{span}\, L)$, then we cannot have $T_{C \cap H}(y) = \mathrm{span}\, L$. $\qquad \square$

*Proof of Theorem 10.* The polar of a convex cone $K$ is

$$K^* = \{\, \delta : \langle w, \delta \rangle \leq 0, \ w \in K \,\}$$

(14, Section 6.E). The double polar theorem (14, Corollary 6.2.1) says that $K^{**} = \mathrm{cl}\, K$. When $K$ is closed, in particular when $K$ is polyhedral, then $K^{**} = K$. Here let $K = \mathrm{con}(\mathrm{pos}(V_{\mathrm{sub}} \setminus \{w\}))$. Then the feasible region for the linear program (19) is $-K^*$. Now the optimal value to (19) is nonpositive if and only if $\langle w, \delta \rangle \leq 0$ for all $\delta \in -K^*$, which is equivalent by the double polar theorem to $w \in (-K^*)^* = -K$ or to $-w \in K$.

Now $w$ is in (18) if and only if $-w$ is a linear combination of elements of $V_{\mathrm{sub}}$ with nonnegative coefficients, that is, if $-w = a \cdot w + \sum_{v \in V_{\mathrm{sub}} \setminus \{w\}} a_v \cdot v$ where $a$ and all the $a_v$ are nonnegative scalars. But this happens if and only if $-w = \sum_{v \in V_{\mathrm{sub}} \setminus \{w\}} (a_v / (1 + a)) \cdot v$, which is equivalent to $-w \in K$. $\qquad \square$

*Proof of Theorem 11.* With probability one

$$Y - y = \sum_{v \in V_{\mathrm{sat}}} b_v(Y) \cdot v$$

where all the coefficients $b_v(Y)$ are nonnegative. From (20a) and (20b) we can derive

$$\langle v, M\delta \rangle = 0, \qquad v \in L_{\mathrm{sat}}$$
$$\langle v, M\delta \rangle < 0, \qquad v \in V_{\mathrm{sat}} \setminus L_{\mathrm{sat}}$$

Hence

$$\langle Y - y, M\delta \rangle = \sum_{v \in V_{\mathrm{sat}} \setminus L_{\mathrm{sat}}} b_v(Y) \cdot \langle v, M\delta \rangle \tag{28}$$

and since all of the $\langle v, M\delta \rangle$ in (28) are strictly negative, the sum can only be zero if all the $b_v(Y)$, $v \in V_{\mathrm{sat}} \setminus L_{\mathrm{sat}}$ are zero. Thus the support of the limiting conditional model consists of points of the form $y + \sum_{v \in L_{\mathrm{sat}}} b_v \cdot v$, where the coefficients are arbitrary. Since all such points are in the preimage of $H_{\mathrm{sub}}$ under the map $y \mapsto M^T y$, we conclude (21) holds. $\qquad \square$